\documentclass[a4paper,12pt]{article}
\usepackage{amssymb}
\usepackage{amsmath}
\usepackage{amstext}

\newtheorem{theorem}{Theorem}

\newtheorem{definition}[theorem]{Definition}

\newtheorem{remark}[theorem]{Remark}

\input tcilatex
\begin{document}

\title{Ultrafunctions and generalized solutions}
\author{Vieri Benci\thanks{
Dipartimento di Matematica Applicata, Universit\`a degli Studi di Pisa, Via
F. Buonarroti 1/c, Pisa, ITALY and Department of Mathematics, College of
Science, King Saud University, Riyadh, 11451, SAUDI ARABIA. e-mail: \texttt{%
benci@dma.unipi.it}}}
\maketitle

\begin{abstract}
The theory of distributions provides generalized solutions for problems
which do not have a classical solution. However, there are problems which do
not have solutions, not even in the space of distributions. As model problem
you may think of%
\begin{equation*}
-\triangle u=u^{p-1}\ ,\ u>0,\ p\geq \frac{2N}{N-2}
\end{equation*}%
with Dirichlet boundary conditions in a bounded open star-shaped set. Having
this problem in mind, we construct a new class of functions called \textbf{%
ultrafunctions} in which the above problem has a (generalized) solution. In
this construction, we apply the general ideas of Non Archimedean Mathematics
(NAM) and some techniques of Non Standard Analysis. Also, some possible
applications of ultrafunctions are discussed.

\medskip

\medskip


\medskip

\end{abstract}

\tableofcontents

\section{Introduction}

We believe that Non Archimedean Mathematics (NAM), namely, mathematics based
on Non Archimedean Fields is very interesting, very rich and, in many
circumstances, allows to construct models of the physical world in a more
elegant and simple way. In the years around 1900, NAM was investigated by
prominent mathematicians such as David Hilbert and Tullio Levi-Civita, but
then it has been forgotten until the '60s when Abraham Robinson presented
his Non Standard Analysis (NSA). We refer to Ehrlich \cite{el06} for a
historical analysis of these facts and to Keisler \cite{keisler76} for a
very clear exposition of NSA.

In this paper we apply the general ideas of NAM and some of the techniques
of NSA to a new notion of generalized functions which we have called \textbf{%
ultrafunctions}. Ultrafunctions are a particular class of functions based on
a superreal field $\mathbb{R}^{\ast }\supset \mathbb{R}$. More exactly, to
any continuous function $f:\mathbb{R}^{N}\rightarrow \mathbb{R}$, we
associate in a canonical way an ultrafunction $f_{\Phi }:\left( \mathbb{R}%
^{\ast }\right) ^{N}\rightarrow \mathbb{R}^{\ast }$ which extends $f;$ but
the ultrafunctions are much more than the functions and among them we can
find solutions of functional equations which do not have any solutions among
the real functions or the distributions.

Now we itemize some of the peculiar properties of the ultrafunctions:

\begin{itemize}
\item the space of ultrafunctions is larger than the space of distributions,
namely, to every distribution $T,$ we can associate in a canonical way an
ultrafunction $T_{\Phi }$ (cf. section \ref{UD});

\item similarly to the distributions, the ultrafunctions are motivated by
the need of having generalized solutions; however, while the distributions
are no longer functions, the ultrafunctions are still functions even if they
have larger domain and range;

\item unlikely the distributions, the space of ultrafunctions is suitable
for non linear problem; in fact any operator $F$ defined for a reasonable
class of functions, can be extended to the ultrafunctions; for example, in
the framework of ultrafunctions $\delta ^{2}\ $makes sense (here $\delta $
is the Dirac measure seen as an ultrafunction);

\item if a problem has a unique classical solution $u,$ then $u_{\Phi }$ is
the only solution in the space of ultrafunctions,

\item the main strategy to prove the existence of generalized solutions in
the space of ultrafunction is relatively simple; it is just a variant of the
Faedo-Galerkin method.
\end{itemize}

This paper is organized as follows. In Section \ref{OT} we introduce NAM via
the notion of $\Lambda $-limit. This approach is quite different from the
usual approach to NAM via NSA. It follows a line developed in \cite{benci95}%
, \cite{benci99}, \cite{BDN2003} and \cite{BHW}. In this section, we
introduce all the notions necessary to understand the rest of the paper, but
we omit details and most of the proofs. In sections \ref{dvb} and \ref{u},
we introduce the notion of ultrafunction and the last three sections are
devoted to applications. The applications are chosen as examples to show the
potentiality of the theory and possible directions of study; they are not an
exhaustive study of the topics treated there.

\bigskip

Before ending the introduction, we want to emphasize the differences by our
approach to NAM and the approach of most people working in Nonstandard
Analysis: there are two main differences, one in the aims and one in the
methods.

Let examine the difference in the aims. We think that infinitesimal and
infinite numbers should not be considered just as entities living in a
parallel universe (the nonstandard universe) which are only a tool to prove
some statement relative to our universe (the standard universe), but rather
that they should be considered mathematical entities which have the same
status of the others and can be used to build models as any other
mathematical entity. Actually, the advantages of a theory which includes
infinitesimals rely more on the possibility of making new models rather than
in the proving techniques. Our papers \cite{BGG} and \cite{BHW} as well as
this one, are inspired by this principle.

As far as the methods are concerned we introduce a non-Archimedean field via
a new notion of limit (see section \ref{OL}). Moreover, we make a very
limited use of logic: the transfer principle (or Leibnitz Principle) is
given by Th. \ref{limit} and it is not necessary to introduce a formal
language. We think that this approach is closer to the way of thinking of
the applied mathematician.

\subsection{Notation}

Let $\Omega $\ be a subset of $\mathbb{R}^{N}$: then

\begin{itemize}
\item $\mathcal{F}\left( \Omega ,E\right) $ denotes the set all the
functions defined in $\Omega $ with values in $E;$

\item $\mathcal{C}\left( \Omega \right) $ denotes the set of real continuous
functions defined on $\Omega ;$

\item $\mathcal{C}_{0}\left( \overline{\Omega }\right) $ denotes the set of
real continuous functions on $\overline{\Omega }$ which vanish on $\partial
\Omega ;$

\item $\mathcal{C}^{k}\left( \Omega \right) $ denotes the set of functions
defined on $\Omega \subset \mathbb{R}^{N}$ which have continuous derivatives
up to the order $k;$

\item $\mathcal{C}_{0}^{k}\left( \overline{\Omega }\right) =\mathcal{C}%
^{k}\left( \overline{\Omega }\right) \cap \mathcal{C}_{0}\left( \overline{%
\Omega }\right) ;$

\item $\mathcal{D}\left( \Omega \right) $ denotes the set of the infinitely
differentiable functions with compact support defined on $\Omega \subset 
\mathbb{R}^{N};\ \mathcal{D}^{\prime }\left( \Omega \right) $ denotes the
topological dual of $\mathcal{D}\left( \Omega \right) $, namely the set of
distributions on $\Omega ;$

\item $\mathcal{S}\left( \Omega \right) $ denotes the Schwartz space and $%
\mathcal{S}^{\prime }\left( \Omega \right) $ the set of tempered
distributions$;$

\item $\mathcal{E}\left( \Omega \right) =\mathcal{C}^{\infty }\left( \Omega
\right) $ denotes the set of the infinitely differentiable functions;\ $%
\mathcal{E}^{\prime }\left( \Omega \right) $ denotes the topological dual of 
$\mathcal{E}\left( \Omega \right) $, namely the set of distributions with
compact support in $\Omega ;$

\item $H^{1}(\Omega )$ is the usual Sobolev space defined as the set of
functions $u\in L^{2}\left( \Omega \right) $ such that $\nabla u\in
L^{2}\left( \Omega \right) ;$

\item $H_{0}^{1}(\Omega )$ is the closure of $\mathcal{D}\left( \Omega
\right) $ in $H^{1}(\Omega );$

\item $H^{-1}(\Omega )$ is the topological dual of $H_{0}^{1}(\Omega ).$
\end{itemize}

\section{$\Lambda $-theory\label{OT}}

As we have already remarked in the introduction, $\Lambda $-theory can be
considered as a variant of nonstandard analysis. It can be introduced via
the notion of $\Lambda $-limit, and it can be easily used for the problems
which we will consider in this paper.

\subsection{Non Archimedean Fields}

In this section, we will give the basic definitions relative to
non-Archimedean fields and some of the basic facts. $\mathbb{F}$ will denote
an ordered field. The elements of $\mathbb{F}$ will be called numbers.
Clearly $\mathbb{F}$ contains (a set isomorphic to) the rational numbers.

\begin{definition}
Let $\mathbb{F}$ be an ordered field. Let $\xi \in \mathbb{F}$. We say that:

\begin{itemize}
\item $\xi $ is infinitesimal if for all $n\in \mathbb{N}$, $|\xi |\ <\frac{1%
}{n}$;

\item $\xi $ is finite if there exists $n\in \mathbb{N}$ such as $|\xi |<n$;

\item $\xi $ is infinite if, for all $n\in \mathbb{N}$, $|\xi |>n$
(equivalently, if $\xi $ is not finite).
\end{itemize}
\end{definition}

\begin{definition}
An ordered field $\mathbb{K}$ is called non-Archimedean if it contains an
infinitesimal $\xi \neq 0$.
\end{definition}

It's easily seen that the inverse of a nonzero infinitesimal number is
infinite, and the inverse of an infinite number is infinitesimal. Clearly,
all infinitesimal numbers are finite.

\begin{definition}
A superreal field is an ordered field $\mathbb{K}$ that properly extends $%
\mathbb{R}$.
\end{definition}

It is easy to show that any superreal field contains infinitesimal and
infinite numbers. Thanks to infinitesimal numbers, in the superreal fields,
we can formalize a new notion of \textquotedblleft closeness".

\begin{definition}
\label{def infinite closeness} We say that two numbers $\xi $ and $\zeta \in 
\mathbb{K}$ are infinitely close if $\xi -\zeta $ is infinitesimal. In this
case, we will write $\xi \sim \zeta $.
\end{definition}

It is easy to see that the relation "$\sim $" of infinite closeness is an
equivalence relation.

\begin{theorem}
If $\mathbb{K}$ is a superreal field, every finite number $\xi \in \mathbb{K}
$ is infinitely close to a unique real number $r\sim \xi $, called the 
\textbf{shadow} or the \textbf{standard part} of $\xi $. We will write $%
r=sh(\xi )$. If $\xi \in \mathbb{K}$ is a positive (negative) infinite
number, then we put $sh(\xi )=+\infty $ ($sh(\xi )=-\infty $).
\end{theorem}

We can also consider the relation of \textquotedblleft finite closeness": 
\begin{equation*}
\xi \sim _{f}\zeta \mathrm{\ if\ and\ only\ if\ }\xi -\zeta \mathrm{\ is\
finite.}
\end{equation*}%
It is readily seen that also $\sim _{f}$ is an equivalence relation. In the
literature, the equivalence classes relative to the two relations of
closeness $\sim $ and $\sim _{f}$, are called monads and galaxies,
respectively.

\begin{definition}
\label{def monad} The monad of a number $\xi $ is the set of all numbers
that are infinitely close to it:%
\begin{equation*}
\mathfrak{mon}(\xi )=\{\zeta \in \mathbb{K}:\xi \sim \zeta \}
\end{equation*}%
The galaxy of a number $\xi $ is the set of all numbers that are finitely
close to it: 
\begin{equation*}
\mathfrak{gal}(\xi )=\{\zeta \in \mathbb{K}:\xi \sim _{f}\zeta \}
\end{equation*}
\end{definition}

So, $\mathfrak{mon}(0)$ is the set of all infinitesimal numbers in $\mathbb{K%
}$ and $\mathfrak{gal}(0)$ is the set of all finite numbers.

\subsection{The $\Lambda $-limit\label{OL}}

$\mathcal{U}$ will denote our "mathematical universe". For our applications
a good choice of $\mathcal{U}$ is given by the superstructure on $\mathbb{R}$%
:

\begin{equation*}
\mathcal{U}=\dbigcup_{n=0}^{\infty }\mathcal{U}_{n}
\end{equation*}%
where $\mathcal{U}_{n}$ is defined by induction as follows:%
\begin{eqnarray*}
\mathcal{U}_{0} &=&\mathbb{R} \\
\mathcal{U}_{n+1} &=&\mathcal{U}_{n}\cup \mathcal{P}\left( \mathcal{U}%
_{n}\right)
\end{eqnarray*}%
Here $\mathcal{P}\left( E\right) $ denotes the power set of $E.$ If we
identify the couples with the Kuratowski pairs and the functions and the
relations with their graphs, clearly $\mathcal{U}$ contains almost all the
mathematical objects needed in mathematics.

Given the universe $\mathcal{U}$, we denote by $\Lambda $ the family of
finite subsets of $\mathcal{U}.$ Clearly $\left( \Lambda ,\subset \right) $
is a directed set and, as usual, a function $\varphi :\Lambda \rightarrow E$
will be called \textit{net }(with values in $E$).

\bigskip {\Large Axioms of\ the }$\Lambda ${\Large -limit}

\begin{itemize}
\item \textsf{(}$\Lambda $-\textsf{1)}\ \textbf{Existence Axiom.}\ \textit{%
There is a superreal field} $\mathbb{K}\supset \mathbb{R}$ \textit{such that
for every net }$\varphi :\Lambda \rightarrow \mathbb{R}$\textit{\ there
exists a unique element }$L\in \mathbb{K\ }$\textit{called the}
\textquotedblleft $\Lambda $-limit" \textit{of}\emph{\ }$\varphi .$ \textit{%
The} $\Lambda $-\textit{limit will be denoted by} 
\begin{equation*}
L=\lim_{\lambda \uparrow \mathcal{U}}\varphi (\lambda )\ \ \text{\textit{or}}%
\ \ L=\lim_{\lambda \in \Lambda }\varphi (\lambda )
\end{equation*}%
\textit{Moreover we assume that every}\emph{\ }$\xi \in \mathbb{K}$\textit{\
is the }$\Lambda $-limit\textit{\ of some net}\emph{\ }$\varphi :\Lambda
\rightarrow \mathbb{R}$\emph{. }

\item ($\Lambda $-2)\ \textbf{Real numbers axiom}. \textit{If }$\varphi
(\lambda )$\textit{\ is} \textit{eventually} \textit{constant}, \textit{%
namely} $\exists \lambda _{0}\in \Lambda :\ \forall \lambda \supset \lambda
_{0},\ \varphi (\lambda )=r,$ \textit{then}%
\begin{equation*}
\lim_{\lambda \uparrow \mathcal{U}}\varphi (\lambda )=r
\end{equation*}

\item ($\Lambda $-3)\ \textbf{Sum and product Axiom}.\ \textit{For all }$%
\varphi ,\psi :\Lambda \rightarrow \mathbb{R}$\emph{: }%
\begin{eqnarray*}
\lim_{\lambda \uparrow \mathcal{U}}\varphi (\lambda )+\lim_{\lambda \uparrow 
\mathcal{U}}\psi (\lambda ) &=&\lim_{\lambda \uparrow \mathcal{U}}\left(
\varphi (\lambda )+\psi (\lambda )\right) \\
\lim_{\lambda \uparrow \mathcal{U}}\varphi (\lambda )\cdot \lim_{\lambda
\uparrow \mathcal{U}}\psi (\lambda ) &=&\lim_{\lambda \uparrow \mathcal{U}%
}\left( \varphi (\lambda )\cdot \psi (\lambda )\right)
\end{eqnarray*}
\end{itemize}

\bigskip

\begin{theorem}
The axioms ($\Lambda $-1)\textsf{,}($\Lambda $-2),($\Lambda $-3) are
consistent.
\end{theorem}

\textbf{Proof. }In order to prove the consistency of these axioms, it is
sufficient to construct a model. Let us consider the algebra $\mathcal{F}%
\left( \Lambda ,\mathbb{R}\right) $ of the real functions defined on $%
\Lambda $ and set%
\begin{equation*}
\mathfrak{I}_{0}=\left\{ \varphi \in \mathcal{F}\left( \Lambda ,\mathbb{R}%
\right) \ |\ \varphi (\lambda )\ \text{is\ eventually}\ 0\right\}
\end{equation*}%
It is easy to check that $\mathfrak{I}_{0}$ is an ideal in the algebra $%
\Lambda .$ By the Krull-Zorn Theorem, every ideal is contained in a maximal
ideal. Let $\mathfrak{I}$ be a maximal ideal containing $\mathfrak{I}_{0}.$
We set 
\begin{equation*}
\mathbb{K}:=\frac{\mathcal{F}\left( \Lambda ,\mathbb{R}\right) }{\cong _{%
\mathfrak{I}}}\ \ \text{\ }
\end{equation*}%
where the equivalence relation $\cong _{\mathfrak{I}}$ is defined as follows:%
\begin{equation*}
\varphi \cong _{\mathfrak{I}}\psi :\Leftrightarrow \varphi -\psi \in 
\mathfrak{I}
\end{equation*}

It is easy to check that $\mathbb{K}$ is an ordered field and $\mathbb{%
R\subset K}$ if we identify $r\in \mathbb{R}$ with the equivalence class $%
\left[ r\right] _{\cong _{\mathfrak{I}}}.$ Finally, we can define the $%
\Lambda $-limit as%
\begin{equation*}
\lim_{\lambda \uparrow \mathcal{U}}\varphi (\lambda )=\left[ \varphi \right]
_{\cong _{\mathfrak{I}}}
\end{equation*}%
Now, it is immediate to check that the $\Lambda $-limit satisfies ($\Lambda $%
-1),($\Lambda $-2),($\Lambda $-3)

$\square $

\bigskip

Now we want to define the $\Lambda $-limit of any bounded net of
mathematical objects in $\mathcal{U}$ (a net $\varphi :\Lambda \rightarrow 
\mathcal{U}$ is called bounded if there exists $n$ such that $\forall
\lambda \in \Lambda ,\varphi (\lambda )\in \mathcal{U}_{n}$). To do this,
consider a net%
\begin{equation}
\varphi :\Lambda \rightarrow \mathcal{U}_{n}  \label{net}
\end{equation}%
We will define $\lim_{\lambda \uparrow \mathcal{U}}\varphi (\lambda )$ by
induction on $n$. For $n=0,$ $\lim_{\lambda \uparrow \mathcal{U}}\varphi
(\lambda )$ is defined by the axioms \textsf{(}$\Lambda $-\textsf{1),}($%
\Lambda $-2),($\Lambda $-3); so by induction we may assume that the limit is
defined for $n-1$ and we define it for the net (\ref{net}) as follows:%
\begin{equation*}
\lim_{\lambda \uparrow \mathcal{U}}\varphi (\lambda )=\left\{ \lim_{\lambda
\uparrow \mathcal{U}}\psi (\lambda )\ |\ \psi :\Lambda \rightarrow \mathcal{U%
}_{n-1},\ \forall \lambda \in \Lambda ,\ \psi (\lambda )\in \varphi (\lambda
)\right\}
\end{equation*}

\begin{definition}
A mathematical entity (number, set, function or relation) which is the $%
\Lambda $-limit of a net is called \textbf{internal}.
\end{definition}

\bigskip

If $E\in \mathcal{U}$, and$\ \varphi :\Lambda \cap \mathcal{P}\left(
E\right) \rightarrow \mathcal{U}_{n},$ then we will use the following
notation:%
\begin{equation*}
\lim_{\lambda \uparrow E}\varphi (\lambda )=\lim_{\mu \uparrow \mathcal{U}%
}\varphi (\mu \cap E).
\end{equation*}%
\bigskip

\subsection{Natural extensions of sets and functions}

\begin{definition}
The \textbf{natural extension }of a set $E\subset \mathbb{R}$ is given by%
\begin{equation*}
E^{\ast }:=\lim_{\lambda \uparrow \mathcal{U}}c_{E}(\lambda )=\ \left\{
\lim_{\lambda \uparrow \mathcal{U}}\psi (\lambda )\ |\ \psi (\lambda )\in
E\right\}
\end{equation*}%
where $c_{E}(\lambda )$ is the net identically equal to $E$.
\end{definition}

Using the above definition we have that 
\begin{equation*}
\mathbb{K}=\mathbb{R}^{\ast }
\end{equation*}

In this context a function $f$ can be identified with its graph; then the
natural extension of a function is well defined. Moreover we have the
following result:

\begin{theorem}
The \textbf{natural extension} of a function%
\begin{equation*}
f:E\rightarrow F
\end{equation*}%
is a function 
\begin{equation*}
f^{\ast }:E^{\ast }\rightarrow F^{\ast };
\end{equation*}%
moreover for every $\varphi :\Lambda \cap \mathcal{P}\left( E\right)
\rightarrow E,$ we have that%
\begin{equation*}
\lim_{\lambda \uparrow \mathcal{U}}\ f(\varphi (\lambda ))=f^{\ast }\left(
\lim_{\lambda \uparrow \mathcal{U}}\varphi (\lambda )\right) .
\end{equation*}
\end{theorem}

When dealing with functions, when the domain of the function is clear from
the context, sometimes the "$\ast $" will be omitted. For example, if $\eta
\in \mathbb{R}^{\ast }$ is an infinitesimal, then clearly $e^{\eta }$ is a
short way to write $\exp ^{\ast }(\eta ).$

The following theorem is a fundamental tool in using the $\Lambda $-limit:

\begin{theorem}
\label{limit}\textbf{(Leibnitz Principle)} Let $\mathcal{R}$ be a relation
in $\mathcal{U}_{n}$ for some $n\geq 0$ and let $\varphi $,$\psi \in 
\mathcal{F}\left( \Lambda ,\mathcal{U}_{n}\right) $. If 
\begin{equation*}
\forall \lambda \in \Lambda ,\ \varphi (\lambda )\mathcal{R}\psi (\lambda )
\end{equation*}%
then%
\begin{equation*}
\left( \underset{\lambda \uparrow \mathcal{U}}{\lim }\varphi (\lambda
)\right) \mathcal{R}^{\ast }\left( \underset{\lambda \uparrow \mathcal{U}}{%
\lim }\psi (\lambda )\right)
\end{equation*}
\end{theorem}

\bigskip

\begin{remark}
Notice that, in the above theorem, the relations "$=$" and "$\in $" do not
change their "meaning", namely "$=^{\ast }$" and "$\in ^{\ast }$" have the
same interpretation than "$=$" and "$\in $".
\end{remark}

\begin{definition}
An internal set is called \textbf{hyperfinite} if it is the $\Lambda $-limit
of finite sets.
\end{definition}

All the internal finite sets are hyperfinite, but there are hyperfinite sets
which are not finite. For example the set%
\begin{equation*}
\mathbb{R}^{\circ }:=\ \underset{\lambda \uparrow \mathcal{U}}{\lim }(%
\mathbb{R}\cap \lambda )
\end{equation*}%
is not finite. The hyperfinite sets are very important since they inherit
many properties of finite sets via Th. \ref{limit}. For example, $\mathbb{R}%
^{\circ }$ has the maximum and the minimum and every internal function%
\begin{equation*}
f:\mathbb{R}^{\circ }\rightarrow \mathbb{R}^{\ast }
\end{equation*}%
has the maximum and the minimum as well.

Also, it is possible to add the elements of an hyperfinite set of numbers or
vectors. Let%
\begin{equation*}
A:=\ \underset{\lambda \uparrow \mathcal{U}}{\lim }A_{\lambda }
\end{equation*}%
be an hyperfinite set; then, the hyperfinite sum is defined as follows: 
\begin{equation*}
\sum_{a\in A}a=\ \underset{\lambda \uparrow \mathcal{U}}{\lim }\sum_{a\in
A_{\lambda }}a
\end{equation*}%
In particular, if $A_{\lambda }=\left\{ a_{1}(\lambda ),...,a_{\beta
(\lambda )}(\lambda )\right\} \ $with\ \ $\beta (\lambda )\in \mathbb{N},\ $%
then, setting 
\begin{equation*}
\beta =\ \underset{\lambda \uparrow \mathcal{U}}{\lim }\ \beta (\lambda )\in 
\mathbb{N}^{\ast }
\end{equation*}
we use the notation%
\begin{equation*}
\sum_{j=1}^{\beta }a_{j}=\ \underset{\lambda \uparrow \mathcal{U}}{\lim }%
\sum_{j=1}^{\beta (\lambda )}a_{j}(\lambda ).
\end{equation*}

\subsection{Qualified sets}

Also, if $Q\subset \Lambda $ and $\varphi :\Lambda \rightarrow \mathcal{U}%
_{n}$, the following notation is quite useful%
\begin{equation*}
\lim_{\lambda \in Q}\varphi (\lambda )=\lim_{\lambda \uparrow \mathcal{U}}%
\widetilde{\varphi }(\lambda )
\end{equation*}%
where 
\begin{equation*}
\widetilde{\varphi }(\lambda )=\left\{ 
\begin{array}{cc}
\varphi (\lambda ) & \text{for}\ \ \lambda \in Q \\ 
\varnothing & \text{for\ }\ \lambda \notin Q%
\end{array}%
\right.
\end{equation*}%
We use this notation to introduce the notion of qualified set:

\begin{definition}
\label{qua}We say that a set $Q\subset \Lambda $ is qualified if for every
bounded net $\varphi ,$ we have that 
\begin{equation*}
\lim_{\lambda \uparrow \mathcal{U}}\varphi (\lambda )=\lim_{\lambda \in
Q}\varphi (\lambda ).
\end{equation*}
\end{definition}

By the above definition, we have that the $\Lambda $-limit of a net $\varphi 
$ depends only on the values that $\varphi $ takes on a qualified set. It is
easy to see that (nontrivial) qualified sets exist. For example, by ($%
\Lambda $-2), we can deduce that, for every $\lambda _{0}\in \Lambda $ the
set%
\begin{equation*}
Q\left( \lambda _{0}\right) :=\left\{ \lambda \in \Lambda \ |\ \lambda
_{0}\subseteq \lambda \right\}
\end{equation*}%
is qualified. In this paper, we will use the notion of qualified set via
this Theorem

\begin{theorem}
\label{billo}Let $\mathcal{R}$ be a relation in $\mathcal{U}_{n}$ for some $%
n\geq 0$ and let $\varphi $, $\psi \in \mathcal{F}\left( \Lambda ,\mathcal{U}%
_{n}\right) $. Then the following statements are equivalent:

\begin{itemize}
\item there exists a qualified set $Q$ such that 
\begin{equation*}
\forall \lambda \in Q,\ \varphi (\lambda )\mathcal{R}\psi (\lambda )
\end{equation*}

\item we have 
\begin{equation*}
\left( \underset{\lambda \uparrow \mathcal{U}}{\lim }\varphi (\lambda
)\right) \mathcal{R}^{\ast }\left( \underset{\lambda \uparrow \mathcal{U}}{%
\lim }\psi (\lambda )\right)
\end{equation*}
\end{itemize}
\end{theorem}

\textbf{Proof}: It is an immediate consequence of Th. \ref{limit} and the
definition of qualified set.

$\square $

\section{The abstract theory\label{dvb}}

In this section we will present a method to extend any vector space $V$ to a
larger vector space $\mathcal{B}\left[ V\right] $ of hyperfinite dimension.
In the next section we will apply this method to functional vector spaces.

\subsection{Definition of ultravectors\label{du}}

\begin{definition}
\label{UV} Let $H\ $be a separable real (or complex) Hilbert space with
scalar product $(\cdot \ ,\cdot )$ and let $V\subset H$ be a dense subspace.
We assume that $H\in \mathcal{U}$ and we set 
\begin{equation*}
\mathcal{B}\left[ V\right] :=\ \underset{\lambda \uparrow V}{\lim }\
V_{\lambda }
\end{equation*}%
where%
\begin{equation*}
V_{\lambda }:=Sp\left( \lambda \right)
\end{equation*}%
is the span of $\lambda $. $\mathcal{B}\left[ V\right] $ is called the space
of ultravectors based on $V.$
\end{definition}

In order to simplify the notation, sometimes, we will set $V_{\mathcal{B}}=%
\mathcal{B}\left[ V\right] .$ Notice that $V_{\mathcal{B}}$ is a vector
space of hyperfinite dimension $\beta \in \mathbb{N}^{\ast }$, were $\beta $
is defined as follows: 
\begin{equation*}
\beta =\dim ^{\ast }(V_{\mathcal{B}})=\ \underset{\lambda \uparrow V}{\lim }%
\left( \dim V_{\lambda }\right) .
\end{equation*}

Let $f\in V;$ if we identify $f$ and $f^{\ast },$ we have that $V\subset V_{%
\mathcal{B}}$. Now let 
\begin{equation}
\Phi :H^{\ast }\rightarrow V_{\mathcal{B}}  \label{prog}
\end{equation}%
be the orthogonal projector. Then, to every vector $f\in H,$ we can
associate the ultravector $\Phi f\in V_{\mathcal{B}}.$ If $\left\{
e_{j}\right\} _{j\leq \beta }$ is a basis for $V_{\mathcal{B}},$then%
\begin{equation}
\Phi f=\sum_{j=1}^{\beta }(f,e_{j})e_{j}  \label{fi}
\end{equation}

Let $V^{\prime }$ denote the dual of $V,$ namely, $V^{\prime }$ is the
family of linear functionals $T$ on $V.$

\begin{definition}
\label{dd}For any $T\in $ $V^{\prime },$ we denote by $\Phi T$ the only
vector in $V_{\mathcal{B}}$ such that 
\begin{equation*}
\forall v\in V_{\mathcal{B}},\ (\Phi T,v)=\left\langle T^{\ast
},v\right\rangle ;
\end{equation*}%
$\Phi T$ is called \textbf{dual} ultravector. Using the orthonormal basis $%
\left\{ e_{j}\right\} _{j\leq \beta }$, we have that%
\begin{equation}
\Phi T=\sum_{j=1}^{\beta }(\Phi T,e_{j})e_{j}=\sum_{j=1}^{\beta
}\left\langle T^{\ast },e_{j}\right\rangle e_{j}  \label{phi}
\end{equation}
\end{definition}

Notice that, if we identify $H$ as a subset of $V^{\prime },$ the operator $%
\Phi $ defined by (\ref{phi}) is the extension of the operator (\ref{fi})
and hence we have denoted them with the same symbol.

From our previous discussion the space of ultravectors $V_{\mathcal{B}}$
contains three types of vectors

\begin{itemize}
\item standard ultravectors: $u\in V_{\mathcal{B}}$ is called \textbf{%
standard} if $u\in V$ (or, to be more precise, if there exists $f\in V$ such
that $u=f^{\ast }$);

\item dual ultravectors: $u\in V_{\mathcal{B}}$ is called \textbf{dual}
ultravector if $u=\Phi T$ for some $T\in V^{\prime }$;

\item proper ultravector: $u\in V_{\mathcal{B}}$ is called \textbf{proper}
ultravector if it is not a dual ultravector.
\end{itemize}

The ultravector which are not standard will be called \textbf{ideal.}

\subsection{Extension of operators\label{EO}}

\begin{definition}
\label{CE}Given the operator $F:D\rightarrow V^{\prime },$ $D\subset V,$ the
map%
\begin{equation*}
F_{\mathbf{\Phi }}:V_{\mathcal{B}}\cap D^{\ast }\rightarrow V_{\mathcal{B}}
\end{equation*}%
defined by%
\begin{equation}
F_{\mathbf{\Phi }}=\Phi \circ F^{\ast }
\end{equation}%
is called \textbf{canonical }extension of $F.$
\end{definition}

By the definition of $F_{\mathbf{\Phi }}$, if $u\in V_{\mathcal{B}}\cap
D^{\ast },$ we have that%
\begin{equation}
\forall v\in V_{\mathcal{B}},\ (F_{\mathbf{\Phi }}\left( u\right)
,v)=\left\langle F^{\ast }\left( u\right) ,v\right\rangle  \label{bellina}
\end{equation}

Using an orthonormal basis $\left\{ e_{j}\right\} _{j\leq \beta }$ for $V_{%
\mathcal{B}},$we have 
\begin{equation*}
F_{\mathbf{\Phi }}\left( u\right) =\sum_{j=1}^{\beta }\left\langle F^{\ast
}(u),e_{j}\right\rangle e_{j}
\end{equation*}

If we identify $H$ with its dual and we take $F:V\cap D\rightarrow H,$ then
equation (\ref{bellina}) becomes:%
\begin{equation}
\forall v\in V_{\mathcal{B}},\ (F_{\mathbf{\Phi }}\left( u\right) ,v)=\left(
F^{\ast }\left( u\right) ,v\right) .  \label{belloccia}
\end{equation}

\section{The ultrafunctions\label{u}}

\bigskip

\subsection{Definition}

\bigskip

\begin{definition}
Let $\Omega $ be a set in $\mathbb{R}^{N}$, and let $V\left( \Omega \right)
\ $be a vector space such that $\mathcal{D}(\Omega )\subseteq V\left( \Omega
\right) \subseteq \mathcal{C}(\Omega )\cap L^{2}(\Omega ).$ Then any function%
\begin{equation*}
u\in \mathcal{B}\left[ V\left( \Omega \right) \right]
\end{equation*}%
is called ultrafunction.
\end{definition}

So the ultrafunctions are $\Lambda $-limits of continuous functions in $V_{%
\mathcal{\lambda }}\left( \Omega \right) :=Sp\left( \lambda \cap V\left(
\Omega \right) \right) $ and hence they are internal functions%
\begin{equation*}
u:\Omega ^{\ast }\rightarrow \mathbb{C}^{\ast }.
\end{equation*}

\begin{remark}
If $V\left( \Omega \right) $ is a Sobolev space such as $H^{1}\left( \Omega
\right) ,$ then the elements of $V\left( \Omega \right) $ are not functions,
but equivalence class of functions, so also the elements of $\mathcal{B}%
\left[ V\left( \Omega \right) \right] $ are equivalence class of functions.
In order to avoid this unpleasant fact, in the definition of ultrafunctions,
we have assumed $V\left( \Omega \right) \subset \mathcal{C}(\Omega )$.
Moreover, this choice has also another motivation: as we will see in the
applications, if we approach a problem via the ultrafunctions, we do not
need Sobolev spaces (even if we might need the Sobolev inequalities). In
some sense the ultrafunctions represent an alternative approach to problems
which do not have classical solutions in some $\mathcal{C}^{k}(\overline{%
\Omega }).$
\end{remark}

Since $V_{\mathcal{B}}(\Omega )\subset \left[ L^{2}(\Omega )\right] ^{\ast },
$ $V_{\mathcal{B}}(\Omega )$ can be equipped with the following scalar
product%
\begin{equation*}
\left( u,v\right) =\int_{\Omega }^{\ast }u(x)\overline{v(x)}\ dx.
\end{equation*}%
where $\int_{\Omega }^{\ast }$ is the natural extension of the Lebesgue
integral considered as a functional.

Notice that the Euclidean structure of $V_{\mathcal{B}}(\Omega )$ is the $%
\Lambda $-limit of the Euclidean structure of every $V_{\lambda }(\Omega )$
given by the usual $L^{2}\left( \Omega \right) $ scalar product.

If $f\in \mathcal{C}(\Omega )$ is a function such that,%
\begin{equation}
\forall g\in V(\Omega ),\ \int f(x)g(x)\ dx<+\infty  \label{rin}
\end{equation}%
then it can be identified with an element of $V\left( \Omega \right)
^{\prime }$ and, by Def. \ref{dd}, there is a unique ultrafunction $f_{\Phi
} $ such that $\forall v\in V_{\mathcal{B}}(\Omega ),$%
\begin{equation}
\int^{\ast }f_{\Phi }(x)v(x)\ dx=\int^{\ast }f^{\ast }(x)v(x)\ dx.
\label{luce}
\end{equation}

The map%
\begin{equation}
\Phi :\mathcal{C}(\Omega )\cap V\left( \Omega \right) ^{\prime }\rightarrow
V_{\mathcal{B}}(\Omega )  \label{canone}
\end{equation}%
is called canonical map. Notice that $f_{\Phi }\neq f^{\ast }$ unless $f\in
V(\Omega ).$

Now let us define a new notion which helps to understand the structure of
ultrafunctions:

\begin{definition}
\label{regular}A hyperfinite basis $\left\{ e_{j}\right\} _{j\leq \beta }$
for $V_{\mathcal{B}}\left( \Omega \right) $ is called \textbf{regular} basis
if

\begin{itemize}
\item it is an orthonormal basis,

\item $\left\{ e_{j}\right\} _{j\in \mathbb{N}}$ is an orthonormal Schauder
basis for $L^{2}(\Omega ).$
\end{itemize}
\end{definition}

\bigskip

The following theorem shows that regular bases exist:

\begin{theorem}
\label{uu}Let $\left\{ h_{j}\right\} _{j\in \mathbb{N}}\subset V\left(
\Omega \right) $ be an orthonormal Schauder basis for $L^{2}(\Omega )$ and
let $W$ be the space generated by \textbf{finite} linear combinations of the
elements of $\left\{ h_{j}\right\} _{j\in \mathbb{N}}$ (hence $W$ is a dense
subspace of $V\left( \Omega \right) ).$ Then there exists a regular basis $%
\left\{ e_{j}\right\} _{j\leq \beta }$ for $V_{\mathcal{B}}\left( \Omega
\right) $ such that%
\begin{equation*}
e_{j}=h_{j}\ \ \text{for }j\leq \theta
\end{equation*}%
where%
\begin{equation*}
\theta =\dim ^{\ast }\left( V_{\mathcal{B}}\left( \Omega \right) \cap
W^{\ast }\right) .
\end{equation*}
\end{theorem}

\textbf{Proof}. Let\textbf{\ }$\left[ \left\{ h_{j}\right\} _{j\in \mathbb{N}%
}\right] ^{\ast }=\left\{ h_{j}\right\} _{j\in \mathbb{N}^{\ast }}\subset
V\left( \Omega \right) ^{\ast }$ be an orthonormal Schauder basis for $%
L^{2}(\Omega )^{\ast }\ $and set 
\begin{equation*}
\theta =\max \left\{ k\in \mathbb{N}^{\ast }\ |\ \forall j\leq k,\ h_{j}\in
V_{\mathcal{B}}\left( \Omega \right) \right\}
\end{equation*}%
Since $\left\{ h_{j}\right\} _{j\in \mathbb{N}}\subset V\left( \Omega
\right) $, $\theta $ is an infinite number in $\mathbb{N}^{\ast }.$ Set $%
e_{j}=h_{j}\ $for $j\leq \theta .$ Now, we can take an orthonormal basis $%
\left\{ e_{j}\right\} _{j\leq \beta }$ for $V_{\mathcal{B}}$ which contains $%
\left\{ e_{j}\right\} _{j\leq \theta }.$

$\square $

So every ultrafunction $u\in V_{\mathcal{B}}\left( \Omega \right) $ can be
represented as follows:%
\begin{equation}
u(x)=\dsum\limits_{j=1}^{\beta }u_{j}e_{j}(x)=\dsum\limits_{n=1}^{\theta
}u_{j}h_{j}(x)+\dsum\limits_{j=\theta +1}^{\beta }u_{j}e_{j}(x)  \label{deco}
\end{equation}%
with%
\begin{equation*}
u_{j}=\int^{\ast }u^{\ast }(x)\overline{e_{j}(x)}\ dx\in \mathbb{R}^{\mathbb{%
\ast }},\ j\leq \beta .
\end{equation*}

In particular, if $f\in L^{2}(\Omega )$ (or more in general if $f\in
V^{\prime }(\Omega )$), the numbers $f_{j},$ $j\in \mathbb{N}$,\ are complex
numbers. The internal function $f_{\Phi }(x)=\dsum\limits_{j=1}^{\beta
}f_{j}e_{j}$ is the orthogonal projection of $f^{\ast }\in L^{2}(\Omega
)^{\ast }$ on $V_{\mathcal{B}}\left( \Omega \right) \subset L^{2}(\Omega
)^{\ast }.$

\bigskip

\textbf{Example}: Let us see an example; we set

\begin{itemize}
\item $\Omega =\left[ 0,1\right] ;$

\item $V\left( \left[ 0,1\right] \right) =C_{0}^{2}\left( \left[ 0,1\right]
\right) ;$

\item $h_{j}(x)=\sqrt{2}\sin \left( j\pi x\right) $;
\end{itemize}

By Th. \ref{uu} there exists a regular basis $\left\{ e_{j}(x)\right\}
_{j\in J}$ which contains $\left\{ \sqrt{2}\sin \left( j\pi x\right)
\right\} _{j\in \mathbb{N}}$. With this assumptions, every vector $u\in V_{%
\mathcal{B}}\left( \left[ 0,1\right] \right) $ can be written as follows%
\begin{equation*}
u(x)=\sqrt{2}\dsum\limits_{n=1}^{\theta }u_{j}\sin \left( j\pi x\right)
+\dsum\limits_{j=\theta +1}^{\beta }u_{j}e_{j}(x)\ \ \text{with}\ \
u_{j}=\int_{0}^{1}u(x)e_{j}(x)dx.
\end{equation*}

\subsection{Ultrafunctions and distributions\label{UD}}

First, we will give a definition of the Dirac $\delta $-ultrafunction
concentrated in $q.$

\begin{theorem}
Given a point $q\in \Omega ,$ there exists a unique function $\delta _{q}$
in $V_{\mathcal{B}}(\Omega )$ such that%
\begin{equation}
\forall v\in V_{\mathcal{B}}(\Omega ),\ \int^{\ast }\delta _{q}(x)v(x)\
dx=v(q).  \label{ddel}
\end{equation}%
$\delta _{q}$ will called the Dirac ultrafunction in $V_{\mathcal{B}}(\Omega
)$ concentrated in $q.$ Moreover, we set $\delta =\delta _{0}.$
\end{theorem}

\textbf{Proof.} Let $\left\{ e_{j}\right\} _{j\leq \beta }$ be any
orthonormal basis for $V_{\mathcal{B}}\left( \Omega \right) $ and set%
\begin{equation*}
\delta _{q}(x)=\dsum\limits_{j=1}^{\beta }e_{j}(q)e_{j}(x)
\end{equation*}%
It is easy to check that $\delta _{q}(x)$ has the desired property; in fact%
\begin{eqnarray*}
\int^{\ast }\delta _{q}(x)v(x)\ dx &=&\int^{\ast }\dsum\limits_{j=1}^{\beta
}e_{j}(q)e_{j}(x)v(x)\ dx \\
&=&\dsum\limits_{j=1}^{\beta }\left( \int^{\ast }e_{j}(x)v(x)\ dx\right)
e_{j}(q)=v(q).
\end{eqnarray*}

$\square $

Next let us see how to associate an ultrafunction $T_{\Phi }=\Phi T$ to
every distribution $T\in \mathcal{D}^{\prime }.$ Let $\left\{ h_{j}\right\}
_{j\in \mathbb{N}}\subset \mathcal{D}$ be an orthonormal Schauder basis for $%
L^{2}(\Omega )$; then, there exists an infinite number $\theta $ such that $%
\left\{ h_{j}\right\} _{j\leq \theta }$ is a basis for $V_{\mathcal{B}%
}(\Omega )\cap \mathcal{D}^{\ast };$ then, $T_{\Phi }(x)$ can be defined as
follows:%
\begin{equation}
T_{\Phi }(x)=\sum_{j=0}^{\theta }\left\langle T^{\ast },h_{j}\right\rangle
h_{j}(x)  \label{expa}
\end{equation}%
Notice that this definition in independent of the choice of the basis since%
\begin{equation}
\int^{\ast }T_{\Phi }(x)\overline{v(x)}\ dx=\left\langle T^{\ast
},v\right\rangle \ \ \text{if}\ \ v\in V_{\mathcal{B}}(\Omega )\cap \mathcal{%
D}^{\ast }  \label{ddd}
\end{equation}%
\begin{equation}
\int^{\ast }T_{\Phi }(x)\overline{v(x)}\ dx=0\ \ \text{if}\ \ v\in \left( V_{%
\mathcal{B}}(\Omega )\cap \mathcal{D}^{\ast }\right) ^{\perp }.  \label{dddd}
\end{equation}%
where $\left( V_{\mathcal{B}}(\Omega )\cap \mathcal{D}^{\ast }\right)
^{\perp }$ denotes the orthogonal complement of $V_{\mathcal{B}}(\Omega
)\cap \mathcal{D}^{\ast }$ in $V_{\mathcal{B}}(\Omega ).$

\begin{remark}
Here the reader must be careful to distinguish the Dirac ultrafunction as
defined by \ref{ddel} and the ultrafunction related to the distribution $%
\delta \in $ $\mathcal{D}^{\prime }$ which now we will call $\delta _{%
\mathcal{D}}.$ In fact, by (\ref{expa}) we have that%
\begin{equation*}
\delta _{\mathcal{D}}(x)=\sum_{j=0}^{\theta }h_{j}(0)h_{j}(x)
\end{equation*}%
while%
\begin{equation*}
\delta (x)=\sum_{j=0}^{\theta }h_{j}(0)h_{j}(x)+\dsum\limits_{j=\theta
+1}^{\beta }e_{j}(0)e_{j}(x)
\end{equation*}%
where $\left\{ h_{j}\right\} _{j\leq \theta }\cup \left\{ e_{j}\right\}
_{\theta +1\leq j\leq \beta }$ is a regular basis for $\left( V_{\mathcal{B}%
}(\Omega )\cap \mathcal{D}^{\ast }\right) ^{\perp }.$ Of course, if $\varphi
\in \mathcal{D}$, we have that%
\begin{equation*}
\int^{\ast }\delta (x)\varphi (x)\ dx=\int^{\ast }\delta _{\mathcal{D}%
}(x)\varphi (x)\ dx=\varphi (0);
\end{equation*}%
actually the above inequality holds for every $\varphi \in V_{\mathcal{B}%
}(\Omega )\cap \mathcal{D}^{\ast }$.
\end{remark}

The above remark suggests the following definition:

\begin{definition}
\label{dt}An ultrafunction $e_{q}\in V_{\mathcal{B}}(\Omega )$ is called a $%
\delta $-type ultrafunction if%
\begin{equation*}
\forall \varphi \in \mathcal{D},\mathcal{\ }\int^{\ast }e_{q}(x)\varphi (x)\
dx\sim \varphi (q).
\end{equation*}
\end{definition}

\bigskip

Following the classification of ultravectors, (\ref{ddd}) and (\ref{dddd}),
the ultrafunctions can be classified as follows:

\begin{definition}
An ultrafunction $u\in V_{\mathcal{B}}(\Omega )$ is called

\begin{itemize}
\item \textbf{standard }if $u\in V(\Omega )$ or, to be more precise, if
there exists $f\in V(\Omega )$ such that $u=f^{\ast }$;

\item \textbf{ideal }if it is not standard;

\item \textbf{dual} ultrafunction if $u=\Phi (T)$ for some $T\in V(\Omega
)^{\prime };$

\item \textbf{distributional ultrafunction} if $u=\Phi (T)$ for some $T\in 
\mathcal{D}^{\prime };$

\item \textbf{proper} ultrafunction if it is not a distributional
ultrafunction.
\end{itemize}
\end{definition}

\section{The Dirichlet problem}

As first application of ultrafunctions, we will consider the following
Dirichlet problem:%
\begin{equation}
\left\{ 
\begin{array}{cc}
u\in \mathcal{C}^{2}(\overline{\Omega }) &  \\ 
-\Delta u=f(x) & \text{for}\ \ x\in \Omega \\ 
u(x)=0 & \text{for\ }x\in \partial \Omega%
\end{array}%
\right.  \label{1}
\end{equation}%
Here $\Omega $ is a bounded set in $\mathbb{R}^{N}.$

This problem is relatively simple and it will help to compare the Sobolev
space approach with the ultrafunctions approach.

\subsection{Generalized solutions\label{gs}}

It is well known that problem (\ref{1}) has a unique solution provided that $%
f(x)\ $and $\partial \Omega $ are smooth. If they are not smooth, it is
necessary to look for generalized solutions. In the Sobolev space approach,
we transform problem (\ref{1}) in the following one:%
\begin{equation}
\left\{ 
\begin{array}{c}
u\in H_{0}^{1}(\Omega ) \\ 
-\Delta u=f(x)%
\end{array}%
\right.  \label{2}
\end{equation}

It is well known that this problem has a unique solution for any bounded
open set $\Omega $ and for a large class of $f,$ namely for every $f\in
H^{-1}(\Omega )$. In this approach, the boundary condition is replaced by
the fact that $u\in H_{0}^{1}(\Omega ),$ namely by the fact that $u$ is the
limit (in $H^{1}(\Omega )$) of a sequence of functions in $\mathcal{C}%
^{2}(\Omega )$ having compact support in $\Omega $. The equation $-\Delta u=f
$ is required to be satisfied in a weak sense:%
\begin{equation*}
-\int_{\Omega }u\Delta \varphi \ dx=\int_{\Omega }f\varphi \ dx\ \ \forall
\varphi \in \mathcal{D}(\Omega )
\end{equation*}%
$u$ itself is not a function but an equivalence class of functions defined $%
a.e.$ in $\Omega .$

\bigskip

Now let us see the ultrafunctions approach. In this case we set $V_{\mathcal{%
B}}^{2,0}(\Omega )=\mathcal{B}\left[ \mathcal{C}_{0}^{2}(\overline{\Omega })%
\right] $ and problem (\ref{1}) can be written as follows:%
\begin{equation}
\left\{ 
\begin{array}{cc}
u\in V_{\mathcal{B}}^{2,0}(\Omega ) &  \\ 
-\Delta _{\Phi }u=f(x) & \text{for}\ \ x\in \text{\ }\Omega ^{\ast }%
\end{array}%
\right.  \label{3}
\end{equation}%
where $\Delta _{\Phi }=\Phi \circ \Delta ^{\ast }:V_{\mathcal{B}%
}^{2,0}(\Omega )\rightarrow V_{\mathcal{B}}^{2,0}(\Omega )$ is given by Def. %
\ref{CE}.

The following result holds:

\begin{theorem}
\label{gugo}For any $f\in V_{\mathcal{B}}^{2,0}(\Omega ),$ problem (\ref{3})
has a unique solution.
\end{theorem}

\textbf{Proof.} By definition, $V_{\mathcal{B}}^{2,0}(\Omega )$ is the $%
\Lambda $-limit of finite dimensional spaces $V_{\lambda }(\Omega )\subset 
\mathcal{C}_{0}^{2}(\overline{\Omega })$. For every $u\in \mathcal{C}%
_{0}^{1}(\overline{\Omega }),$ by the Poincar\'{e} inequality, we have that 
\begin{equation*}
\int_{\Omega }\nabla u\cdot \nabla u\ dx\geq k\left\Vert u\right\Vert
_{L^{2}(\Omega )}^{2}.
\end{equation*}

In particular, the above inequality holds for any $u\in V_{\lambda }(\Omega )
$. Now, let 
\begin{equation*}
\Phi _{\lambda }:L^{2}(\Omega )\rightarrow V_{\lambda }(\Omega ),
\end{equation*}%
be the orthogonal projection. For every $u,v\in V_{\lambda }(\Omega ),$ we
have that%
\begin{equation*}
\int_{\Omega }\nabla u\cdot \nabla v\ dx=\int_{\Omega }-\Delta u\ v\ dx
\end{equation*}%
Then, by the Poincar\'{e} inequality, 
\begin{equation*}
-\Phi _{\lambda }\Delta :V_{\lambda }(\Omega )\rightarrow V_{\lambda
}(\Omega )
\end{equation*}%
is a positive definite symmetric operator. Then it is invertible. So we have
that, for any $\lambda \in \Lambda ,$ there exists a unique $\bar{u}%
_{\lambda }\in V_{\lambda }(\Omega )$ such that%
\begin{equation}
\forall v\in V_{\lambda }(\Omega ),\ \int_{\Omega }-\Delta \bar{u}_{\lambda
}v\ dx=\int_{\Omega }f_{\lambda }v\ dx  \label{robina}
\end{equation}%
where $f_{\lambda }\in V_{\lambda }(\Omega )$ is such that $f=\ \underset{%
\lambda \uparrow V}{\lim }\ f_{\lambda }.$ If we take the $\Lambda $-limit
in this equality, we get%
\begin{equation}
\forall v\in V_{\mathcal{B}}^{2,0}(\Omega ),\ -\int_{\Omega }^{\ast }\Delta
^{\ast }\bar{u}\ v\ dx=\int_{\Omega }^{\ast }fv\ dx  \label{roba}
\end{equation}%
where%
\begin{equation*}
\bar{u}=\ \underset{\lambda \uparrow V}{\lim }\ \bar{u}_{\lambda }
\end{equation*}%
and hence, by (\ref{belloccia}), we get 
\begin{equation*}
-\Delta _{\Phi }\bar{u}=f
\end{equation*}

The uniqueness follows from the uniqueness of $\bar{u}_{\lambda }$.

$\square $

\begin{remark}
This example shows quite well the general strategy to solve problems within
the framework of ultrafunctions. First you solve a finite dimensional
problem and then you take the $\Lambda $-limit. Since the $\Lambda $-limit
exists for any sequence of mathematical objects, the solvability of the
finite dimensional approximations imply the existence of a generalized
solution.
\end{remark}

The solution is a function $\bar{u}:\overline{\Omega }^{\ast }\rightarrow 
\mathbb{R}^{\ast };$ $\bar{u}$ is defined for every $x\in \overline{\Omega }%
^{\ast },$ and we have that $u(x)=0$ for\ $x\in \partial \Omega ^{\ast }.$
So the boundary condition can be interpreted "classically" while this is not
possible in $H_{0}^{1}(\Omega )$. If problem (\ref{1}) has a solution $U\in 
\mathcal{C}^{2}(\overline{\Omega }),$ then 
\begin{equation*}
\bar{u}=U^{\ast }.
\end{equation*}%
If problem (\ref{2}) has a solution $U\in H_{0}^{1}(\Omega ),$ then we have
that 
\begin{equation*}
\int_{\Omega }U\varphi \ dx\sim \int_{\Omega }^{\ast }\bar{u}\varphi \ dx\ \
\forall \varphi \in \mathcal{C}_{0}^{2}(\overline{\Omega })
\end{equation*}

Notice that in the above formula the left hand side integral is a Lebesgue
integral while in the right hand side, $\int^{\ast }$ is the $\ast $%
-transform of the Riemann integral; the integral make sense since $\bar{u}%
,\varphi \in \left[ \mathcal{C}_{0}(\overline{\Omega })\right] ^{\ast }$. In
the theory of ultrafunctions, the Lebesgue integral seems to be not so
necessary.

There are interesting and physically relevant cases in which the
generalization of the Dirichlet problem cannot be treated within the Sobolev
space $H_{0}^{1}(\Omega ).$ For example, consider the problem:%
\begin{equation}
\left\{ 
\begin{array}{cc}
-\Delta u=\delta _{y} & \text{for}\ \ x\in \Omega \\ 
u(x)=0 & \text{for\ }x\in \partial \Omega%
\end{array}%
\right.  \label{4}
\end{equation}%
where $\delta _{y}$ is the Dirac measure concentrated at $y\in \Omega $.
This problem is quite natural in potential theory; in fact $u$ represents
the potential generated by a point source (and usually it is called Green
function). However this problem does not have solution in $H_{0}^{1}(\Omega
) $ since $\delta \notin H^{-1}(\Omega ).$ Actually, with some work, it is
possible to prove that it has a "generalized solution" in $H_{0}^{1}(\Omega
)+\mathcal{E}^{^{\prime }}(\Omega ).$ However, in the framework of
ultrafunction, problem (\ref{4}) is nothing else but a particular case of
problem (\ref{3}).

However, if $f\in V_{\mathcal{B}}^{2,0}(\Omega )$ is a proper ultrafunction,
(namely, $f$ cannot be associated to a distribution via (\ref{ddd}) and (\ref%
{dddd})), problem (\ref{3}) has a solution which cannot be interpreted as a
distribution solution. For example, you can take $f=\delta (x)^{2}.$
Remember that $\delta (x)^{2},$ in the ultrafunction theory, makes sense by
Def. \ref{CE}.

\begin{remark}
If you take $f=\delta ^{2}$ you get a well posed mathematical problem, but,
most likely, it does not represent any "physically" relevant phenomenon.
However, it is possible to choose some proper ultrafunction $f\in V_{%
\mathcal{B}}^{2,0}(\Omega )$ which models physical phenomena. For example%
\begin{equation*}
f(x)=\sin \alpha (\mathbf{n}\cdot x);\ \mathbf{n}\in \mathbb{R}^{N},\
\left\vert \mathbf{n}\right\vert =1,\ \alpha \in \mathbb{R}^{\ast }\ \text{%
infinite, }x\in K^{\ast },K\subset \subset \Omega
\end{equation*}%
might represent a electrostatic problem in a sort of periodic medium such as
a crystal. Here $K$ represent the support of the crystal and $f(x)$
represents its charge density; it consists of periodic layers of positive
and negative charges at a distance of $\frac{1}{\pi \alpha }.$ From a
macroscopic point of view the solution is $0$, but at the microscopic level
this is not the case. In fact the solution $u$ of problem (\ref{3}) does not
vanish, even if it can be proved that 
\begin{equation*}
\forall v\in \mathcal{C}^{2}(\overline{\Omega }),\ \int_{\Omega }^{\ast }%
\bar{u}\ v\ dx\sim 0.
\end{equation*}
\end{remark}

\subsection{The variational approach}

Looking at problem (\ref{1}) from a variational point of view, the
comparison between the Sobolev space approach and the ultrafunctions
approach becomes richer.

\bigskip

It is well known that the equation (\ref{1}) is the Euler-Lagrange equation
of the energy functional%
\begin{equation*}
J(u)=\int_{\Omega }\left( \frac{1}{2}\left\vert \nabla u\right\vert
^{2}-fu\right) \ dx
\end{equation*}%
Thus a minimizer of $J(u)$ on $\mathcal{C}_{0}^{2}(\overline{\Omega })$
solves the problem. However, if $f(x)\ $and $\partial \Omega $ are not
smooth a minimizing sequence does not converge in $\mathcal{C}_{0}^{2}(%
\overline{\Omega })$ and also when it converges, it can be proved only by
making hard estimates.

On the other hand, if you define $H_{0}^{1}(\Omega )$ as the closure of $%
\mathcal{D}(\Omega )$ with respect to the norm 
\begin{equation*}
\left\Vert u\right\Vert _{H_{0}^{1}}=\sqrt{\int_{\Omega }\left\vert \nabla
u\right\vert ^{2}dx}
\end{equation*}%
the functional $\ J(u)$ becomes $\frac{1}{2}\left\Vert u\right\Vert
_{H_{0}^{1}}^{2}-\int_{\Omega }fu\ dx$ and it is immediate to see that it
has a minimizer provided that $f\in H^{-1}(\Omega ).$

If you consider problem (\ref{4}), the trouble with the energy functional is
that the energy%
\begin{equation*}
J(u)=\int_{\Omega }\frac{1}{2}\left\vert \nabla u\right\vert ^{2}dx-u(y),\ \
u\in \mathcal{C}_{0}^{2}(\overline{\Omega })
\end{equation*}%
is not bounded below and $J$ cannot be extended to all $H_{0}^{1}(\Omega ).$

Instead, if we use the ultrafunctions approach, the energy%
\begin{equation*}
J(u)=\int_{\Omega ^{\ast }}^{\ast }\left( \frac{1}{2}\left\vert \nabla
u\right\vert ^{2}-\delta _{y}u\right) dx,\ \ u\in V_{\mathcal{B}%
}^{2,0}(\Omega )
\end{equation*}%
is well defined and it makes sense to look for a minimizer in $V_{\mathcal{B}%
}^{2,0}(\Omega ).$ For every $\lambda \subset \mathcal{C}_{0}^{2}(\overline{%
\Omega })\cap \Lambda ,\ J(u)$ has a minimizer $u_{\lambda }$ in $V_{\lambda
}(\Omega )\subset \mathcal{C}_{0}^{2}(\overline{\Omega })$, and hence, if
you set 
\begin{equation*}
\bar{u}=\ \underset{\lambda \uparrow V}{\lim }\ u_{\lambda },
\end{equation*}%
we have that 
\begin{equation*}
J(\bar{u})=\ \underset{\lambda \uparrow V}{\lim }\left[ \int_{\Omega }\frac{1%
}{2}\left\vert \nabla u_{\lambda }\right\vert ^{2}\ dx-u_{\lambda }(y)\right]
\end{equation*}%
minimizes $J(u)$ in $V_{\mathcal{B}}^{2,0}(\Omega ).$ Clearly, for some
values of $u$, $J(u)$ may assume infinite values in $\mathbb{R}^{\ast }$,
but this is not a problem, actually in my opinion, this is one of the main
reason to legitimate the use non-Archimedean fields. In fact in the
framework of NAM, it is possible to make models of the physical world in
which there are material points with a finite charge. They have an
"infinite" energy, but, nevertheless, we can make computations and if
necessary to evaluate it. The epistemological (and very interesting) issue
relative to the meaning of their "physical existence" should not prevent
their use.

\section{The bubbling phenomenon relative to the Sobolev critical exponent}

The bubbling phenomenon relative to the critical Sobolev exponent is the
model problem which has inspired this work. In general (at least in the
simplest cases), the bubbling phenomenon consists in minimizing sequences
whose \textit{mass concentrate to some points}; however their "limit" does
not exist in any Sobolev space and not even in any distribution space due to
the "strong" non-linearity of the problem. Nevertheless, these problems have
been extensively studied and we know a lot of facts relative to the
minimizing sequences (or more in general to non-converging Palais-Smale
sequences) which, up to an equivalence relation, are called \textit{critical
points at infinity} (see \cite{bahri}). The literature on this topic is huge
(you can find part of it in \cite{Cha}). We refer also to \cite{bahri}, \cite%
{BN} and \cite{Cha} for an exposition of the utility of knowing the
properties of the critical points at infinity.

Ultrafunction theory seems to be an appropriate tool to deal with these kind
of problems.

\subsection{Description of the problem}

Let us consider the following minimization problem:%
\begin{equation*}
\underset{u\in \mathfrak{M}_{p}}{\min }\ J(u)
\end{equation*}%
where%
\begin{equation*}
J(u)=\int_{\Omega }\left\vert \nabla u\right\vert ^{2}\ dx
\end{equation*}%
and 
\begin{equation*}
\mathfrak{M}_{p}=\left\{ u\in \mathcal{C}_{0}^{2}(\overline{\Omega }):\
\int_{\Omega }\left\vert u\right\vert ^{p}\ dx=1\right\}
\end{equation*}%
Here $\Omega $ is a bounded set in $\mathbb{R}^{N}\ $with smooth boundary,\ $%
N\geq 3$ and $p>2$. If $J$ has a minimizer, it is a solution of the
following elliptic eigenvalue problem: 
\begin{equation}
\left\{ 
\begin{array}{cc}
u\in \mathcal{C}_{0}^{2}(\overline{\Omega }) &  \\ 
-\Delta u=\lambda u^{p-1} & \text{for}\ \ x\in \Omega \\ 
u(x)>0 & \text{for}\ \ x\in \Omega \\ 
\int_{\Omega }\left\vert u\right\vert ^{p}dx=1 & 
\end{array}%
\right.  \label{YAM}
\end{equation}

As usual in the literature, we set%
\begin{equation*}
2^{\ast }=\frac{2N}{N-2};
\end{equation*}%
$2^{\ast }$ is called the critical Sobolev exponent for problem (\ref{YAM})
(notice that this "$\ast $" has nothing to do with the natural extension).
Moreover, we set%
\begin{equation*}
m_{p}:=\ \underset{u\in \mathfrak{M}_{p}}{\inf }\ J(u)
\end{equation*}%
The following facts are well known (see e.g. \cite{Cha} and references):

\begin{itemize}
\item (i) if $2<p<2^{\ast },$ then $m_{p}>0$ and it is achieved; hence
problem (\ref{YAM}) has a solution.

\item (ii) if $p=2^{\ast },$ then $m_{2^{\ast }}>0$ and it is achieved only
if $\Omega =\mathbb{R}^{N};$ however there are particular domains $\Omega $
such that (\ref{YAM}) has a solution (which, of course, is not a minimizer
of $J,$ but a critical point).

\item (iii) if $p>2^{\ast },$ then $m_{p}=0$ and it is not achieved.
\end{itemize}

Probably, the most interesting case is the second one (the critical exponent
case) since it presents many interesting phenomena. If $u_{n}$ is a
minimizing sequence, it has a subsequence $u_{n}^{\prime }$ which
concentrates to some point $x_{0}\in \overline{\Omega };$ more exactly, $%
u_{n}^{\prime }\rightharpoonup 0$ weakly in $H_{0}^{1}(\Omega )\ $and
strongly in $H_{0}^{1}(\Omega \setminus B_{\varepsilon }(x_{0}));$
consequently, $\left\vert u_{n}^{\prime }\right\vert ^{p}\rightharpoonup
\delta _{x_{0}}$ weakly in $\mathcal{D}^{\prime }(\Omega ),$ but $\left(
\delta _{x_{0}}\right) ^{1/p}$ cannot be interpreted as a generalized
solution in the framework of the distribution theory just because $\left(
\delta _{x_{0}}\right) ^{1/p}$ makes no sense. This phenomenon is called
"bubbling" and probably problem (\ref{YAM}) with $p=2^{\ast }$ is the
simplest problem which presents it. Similar phenomena occur in many other
variational problems such as the Yamabe problem, the Kazdan-Warner problem,
in the study of harmonic maps between manifolds, in minimal surfaces theory,
in the Yang-Mills equations etc.

Let us go back to discuss the concentration phenomenon of a minimizing
sequence. Not all the points of $\overline{\Omega }$ have the same "dignity"
as concentration points. Let us explain what do we mean.

Let 
\begin{equation}
u_{p},\ p\in (2,2^{\ast }),  \label{grulla}
\end{equation}%
be a minimizer of $J(u)$ on the set $\mathfrak{M}_{p}$. If $p\rightarrow
2^{\ast }$ from the left, it is well known that%
\begin{equation*}
\underset{p\rightarrow \left( 2^{\ast }\right) ^{-}}{\lim }m_{p}=m_{2^{\ast
}}
\end{equation*}%
and that 
\begin{equation}
v_{p}:=\frac{u_{p}}{\int_{\Omega }\left\vert u_{p}\right\vert ^{2^{\ast }}dx}%
,  \label{lilla}
\end{equation}%
is a minimizing sequence of $J$ on $\mathfrak{M}_{2^{\ast }}.$ If, for every 
$u\in \mathfrak{M}_{2^{\ast }},$ we set%
\begin{equation*}
\mathfrak{B}\left( u\right) =\int_{\Omega }x\left\vert u\right\vert
^{2^{\ast }}dx
\end{equation*}%
then we have that, in the generic case, 
\begin{equation*}
\underset{p\rightarrow \left( 2^{\ast }\right) ^{-}}{\lim }\mathfrak{B}%
\left( v_{p}\right) =\overline{x}
\end{equation*}%
where $\overline{x}$ is an interior point of $\Omega $. Thus, in this sense, 
$\overline{x}$ is a "special" concentration point. If we apply ultrafunction
theory, the world "special" will get a new meaning; in fact $\overline{x}$
will be characterized as the point infinitely close to the concentration
point of the generalized solution. This issue will be further discussed in
the next section.

\bigskip

\subsection{Generalized solutions}

\bigskip

The minimization problem considered in the previous section can be studied
in the framework of the ultrafunctions. In this framework the problem takes
the following form:

\begin{equation}
\underset{u\in \widetilde{\mathfrak{M}}_{p}}{\min }\ J(u)  \label{*}
\end{equation}%
where%
\begin{equation*}
J(u)=\int_{\Omega }^{\ast }\left\vert \nabla u\right\vert ^{2}\ dx
\end{equation*}%
and 
\begin{equation*}
\widetilde{\mathfrak{M}}_{p}=\left\{ u\in V_{\mathcal{B}}^{2,0}(\overline{%
\Omega })\ |\ \int_{\Omega }^{\ast }\left\vert u\right\vert ^{p}\
dx=1\right\}
\end{equation*}%
where $V_{\mathcal{B}}^{2,0}(\overline{\Omega })=\mathcal{B}\left[ \mathcal{C%
}_{0}^{2}(\overline{\Omega })\right] .$

\begin{theorem}
For every $p>2,$ problem (\ref{*}) has a solution $\tilde{u}_{p}.$ If we set 
$\widetilde{m}_{p}=J(\tilde{u}_{p}),$ we have the following

\begin{itemize}
\item (i) if $2<p<2^{\ast },$ then $\widetilde{m}_{p}=m_{p}\in \mathbb{R}%
^{+} $ and there is at least one standard minimizer $\tilde{u}_{p}$, namely $%
\tilde{u}_{p}\in \mathcal{C}_{0}^{2}(\overline{\Omega });$

\item (ii) if $p=2^{\ast },$(and $\Omega \neq \mathbb{R}^{N}),$ then $%
\widetilde{m}_{2^{\ast }}=m_{2^{\ast }}+\varepsilon $ where $\varepsilon $
is a positive infinitesimal;

\item (iii) if $p>2^{\ast },$ then $\widetilde{m}_{p}=\varepsilon _{p}$
where $\varepsilon _{p}$ is a positive infinitesimal.
\end{itemize}
\end{theorem}

\bigskip

\textbf{Proof.} The proof of this theorem is a simple application of the
nostandard methods. We will describe it with some details for the reader not
acquainted with these methods.

We set%
\begin{equation*}
\tilde{u}_{p}=\ \underset{\lambda \uparrow \mathcal{C}_{0}^{2}(\overline{%
\Omega })}{\lim }\ u_{p,\lambda }
\end{equation*}%
where $u_{p,\lambda }$ is the minimizer of $J(u)$ on the set $\mathfrak{M}%
_{p}\cap V_{\lambda }(\overline{\Omega });\ V_{\lambda }(\overline{\Omega }%
)=Sp(\lambda )\subset \mathcal{C}_{0}^{2}(\overline{\Omega })$. We recall
that $\mathfrak{M}_{p}\cap V_{\lambda }(\overline{\Omega })\neq 0$ for $%
\lambda $ in a qualified set and that the minimum exists since $V_{\lambda }(%
\overline{\Omega })$ is a finite dimensional vector space and hence $%
\mathfrak{M}_{p}\cap V_{\lambda }(\overline{\Omega })$ is compact. If we set 
\begin{equation*}
m_{p,\lambda }:=\underset{u\in \mathfrak{M}_{p}\cap V_{\lambda }(\overline{%
\Omega })}{\min }\ J(u),
\end{equation*}%
taking the $\Lambda $-limit, we have that 
\begin{equation*}
\widetilde{m}_{p}:=\underset{\lambda \uparrow \mathcal{C}_{0}^{2}(\overline{%
\Omega })}{\lim }\ m_{p,\lambda }=\ \underset{u\in \widetilde{\mathfrak{M}}%
_{p}}{\min }\ J(u).
\end{equation*}%
So the existence result is proved. Now let us prove the second part of the
theorem:

(i) If you take $\lambda _{0}=\left\{ u_{p}\right\} ,$ where $u_{p}$ is
given by (\ref{grulla}) then for every $\lambda \supseteq \lambda _{0}$, we
have that 
\begin{equation*}
m_{p,\lambda }:=\underset{u\in \mathfrak{M}_{p}\cap V_{\lambda }(\overline{%
\Omega })}{\min }J(u)=J(u_{p})=m_{p}
\end{equation*}%
and hence, taking the $\Lambda $-limit, we have that $\widetilde{m}%
_{p}=m_{p} $.

(ii) It is well known that the value $\widetilde{m}_{2^{\ast }}$ is not
achieved by any function $u\in \mathfrak{M}_{2^{\ast }}\cap V_{\lambda }(%
\overline{\Omega });$ then $m_{2^{\ast },\lambda }>m_{2^{\ast }},$ and
hence, taking the $\Lambda $-limit, we have that $\widetilde{m}_{2^{\ast
}}>m_{2^{\ast }}.$ On the other hand, for every $b\in \mathbb{R}^{+},$ there
exists $u\in \mathfrak{M}_{2^{\ast }}$ such $J(u)\leq m_{2^{\ast }}+b,$ and
hence 
\begin{equation*}
\widetilde{m}_{2^{\ast }}=J(\tilde{u}_{2^{\ast }})\leq J(u)\leq m_{2^{\ast
}}+b,
\end{equation*}%
and so, by the arbitrariness of $b,$ we get that $\widetilde{m}_{2^{\ast
}}\sim m_{p}.$

(iii) follows by the same argument used in (ii) replacing $m_{2^{\ast }}$
with $0.$

\bigskip $\square $

The next theorem shows that, for $p=2^{\ast },$ the solution $\tilde{u}$
concentrates where it is expected to do.

\bigskip

\begin{theorem}
Suppose that problem (\ref{YAM}) (with $p=2^{\ast }$) has a unique minimum $%
\tilde{u}$ and set 
\begin{equation*}
\xi =\mathfrak{B}^{\ast }\left( \tilde{u}\right) :=\int_{\Omega }^{\ast
}x\left\vert \tilde{u}\right\vert ^{2^{\ast }}dx\in \Omega ^{\ast }.
\end{equation*}%
Then%
\begin{equation*}
\xi \sim \underset{p\rightarrow \left( 2^{\ast }\right) ^{-}}{\lim }%
\mathfrak{B}\left( v_{p}\right) .
\end{equation*}%
where $v_{p}$ is defined by (\ref{lilla}).
\end{theorem}

\textbf{Proof.} Fix $r\in \mathbb{R}^{+}$. We want to prove that, for $p$
sufficiently close to $2^{\ast },$ we have that%
\begin{equation*}
d^{\ast }(\mathfrak{B}\left( v_{p}\right) ,\xi )\leq r
\end{equation*}%
where $d^{\ast }$ denotes the distance in $\left( \mathbb{R}^{N}\right)
^{\ast }$. We have that 
\begin{equation}
\xi =\ \underset{\lambda \uparrow \mathcal{C}_{0}^{2}(\overline{\Omega })}{%
\lim }\ x_{\lambda }  \label{lulu}
\end{equation}%
where $x_{\lambda }=\mathfrak{B}\left( u_{\lambda }\right) $ and $u_{\lambda
}$ is a minimizer of $J$ on the manifold $\mathfrak{M}_{2^{\ast }}\cap
V_{\lambda }$. Let $\tilde{u}$ be the minimum of $J\ $on $\widetilde{%
\mathfrak{M}}_{2^{\ast }}$, and apply Th.\ \ref{billo}, to the relation $%
\mathcal{R}$ defined as follows:%
\begin{equation*}
u_{\lambda }\mathcal{R}\left( \mathfrak{M}_{2^{\ast }}\cap V_{\lambda
}(\Omega )\right)
\end{equation*}%
if and only if%
\begin{equation*}
u_{\lambda }\text{\ is the unique minimum of }J\text{\ on }\mathfrak{M}%
_{2^{\ast }}\cap V_{\lambda }(\Omega ).
\end{equation*}

Then by Th.\ \ref{billo}, there exists a qualified set $Q\subset \Lambda
(V), $ such that, for every $\lambda \in Q,\ u_{\lambda }$ is the unique
minimum of $J\ $on $\mathfrak{M}_{2^{\ast }}\cap V_{\lambda }(\Omega ).$

Thus $\exists b\in \mathbb{R}^{+},$ $\exists \lambda _{0},\forall \lambda
\geq \lambda _{0},\lambda \in Q,\forall u\in \mathfrak{M}_{2^{\ast }}\cap
V_{\lambda }$%
\begin{equation*}
J(u)<m_{2^{\ast }}+b\Rightarrow d^{\ast }(\mathfrak{B}\left( u\right)
,x_{\lambda })\leq \frac{r}{2}
\end{equation*}%
and hence, may be taking a bigger $\lambda _{0},$ using (\ref{lulu}), we get%
\begin{equation}
J(u)<m_{2^{\ast }}+b\Rightarrow d^{\ast }(\mathfrak{B}\left( u\right) ,\xi
)\leq r  \label{lella2}
\end{equation}

Now, let $v_{p}$ be the function defined by (\ref{lilla}); it is well known
that 
\begin{equation*}
\underset{p\rightarrow \left( 2^{\ast }\right) ^{-}}{\lim }\
J(v_{p})=m_{2^{\ast }}
\end{equation*}%
Then we can take $p$ so close to $2^{\ast }$ so that 
\begin{equation*}
J(v_{p})\leq m_{2^{\ast }}+b.
\end{equation*}%
Since $v_{p}\in \mathfrak{M}_{2^{\ast }}\cap V_{\lambda },$ for every $%
\lambda \geq \lambda _{0}\cup \left\{ v_{p}\right\} ,$ $\lambda \in Q$, by (%
\ref{lella2}), we get that 
\begin{equation*}
d^{\ast }(\mathfrak{B}\left( v_{p}\right) ,\xi )\leq r.
\end{equation*}

$\square $

\bigskip

\begin{remark}
If $J$ does not have a unique minimum, but a set of minimizers, we set%
\begin{equation*}
\Gamma =\left\{ \xi \in \Omega ^{\ast }:\ \xi =\mathfrak{B}\left( \tilde{u}%
\right) \ \text{where }\tilde{u}\ \text{is a minimizer}\right\} .
\end{equation*}
Then, arguing as in the proof of the above theorem, it is easy to get the
following result: let $p_{n}\rightarrow (2^{\ast })^{-}$, let $x_{n}=%
\mathfrak{B}\left( v_{p_{n}}\right) $ and let $x_{n}^{\prime }$ be a
converging subsequence of $x_{n}$. Then there exists $\xi \in \Gamma $ such
that%
\begin{equation*}
\xi \sim \ \underset{n\rightarrow \infty }{\lim }x_{n}^{\prime }
\end{equation*}
\end{remark}

\section{Ultrafunctions and Quantum Mechanics}

In this section we will describe an application of the previous theory to
the formalism of Quantum Mechanics. In the usual formalism, a physical state
is described by a unit vector $\psi \ $in a Hilbert space $\mathcal{H}$ and
an observable by a self-adjoint operator defined on it. In the
ultravectors/ultrafunctions formalism, a physical state is described by a
unit vector $\psi \ $in a hyperfinite space of ultravectors $V_{\mathcal{B}}$
and an observable by a Hermitian operator defined on it.

We think that the ultravectors approach presents the following advantages:

\begin{itemize}
\item once you have learned the basic facts of the $\Lambda $-theory, the
formalism which you get is easier to handle since it is based on the matrix
theory on finite vector spaces rather than on unbounded self-adjoint
operators in Hilbert spaces;

\item this approach is closer to the "infinite" matrix approach of the
beginning of QM before the work of von Neumann and also closer to the way of
thinking of the theoretical physicists and chemists;

\item all observables (hyperfinite matrices) have infinitely many
eigenvectors; so the continuous spectrum can be considered as a set of
eigenvalues infinitely close to each other;

\item the distinction between standard and ideal ultravectors has a physical
meaning;

\item the dynamics does not present any difficulty since it is given by the
exponential matrix relative to the Hamiltonian matrix.
\end{itemize}

Clearly it is too early to know if this formalism will lead to some new
physically relevant fact; in any case we think that it is worthwhile to
investigate it. In this paper we limit ourselves only to some very general
remark.

\subsection{The axioms of Quantum Mechanics}

We start giving a list of the main axioms of quantum mechanics as it is
usually given in any textbook and then we will compare it with the
alternative formalism based on ultravectors.

\bigskip

{\Large Classical axioms of QM}

\bigskip

\textbf{Axiom C1}. A physical state is described by a unit vector $\psi \ $%
in a Hilbert space $\mathcal{H}$.

\bigskip

\textbf{Axiom C2.} An observable is represented by a self-adjoint operator $%
A $ on $\mathcal{H}$.

(a) The set of observable outcomes is given by the eigenvalues $\mu _{j}$ of 
$A$.

(b) After an observation/measurement of an outcome $\mu _{j}$, the system is
left in a eigenstate $\psi _{j}$ associated with the detected eigenvalue $%
\mu _{j}$.

(b) In a measurement the transition probability $\mathcal{P}$ from a state $%
\psi $ to an eigenstate $\psi _{j}$ is given by 
\begin{equation*}
\mathcal{P}=\left\vert \left( \psi ,\psi _{j}\right) \right\vert ^{2}.
\end{equation*}

\textbf{Axiom C3}. The evolution of a state is given by the Shroedinger
equation%
\begin{equation*}
i\frac{\partial \psi }{\partial t}=H\psi
\end{equation*}%
where $H,\ $the Hamiltonian operator, is a self-adjoint operator
representing the energy of the system.

\bigskip

{\Large Axioms of QM based on ultravectors}

\bigskip

\textbf{Axiom U1}. A physical system is described by a complex
valued-ultravector space $V_{\mathcal{B}}=\mathcal{B}\left[ V\right] ;\ $a
state of this system is described by a unit ultravector vector $\psi \ $in $%
V_{\mathcal{B}}$.

\bigskip

\textbf{Axiom U2.} An observable is represented by a Hermitian operator $A$
on $V_{\mathcal{B}}$.

(a) The set of observable outcomes is given by $sh\left( \mu _{j}\right) $
where $\mu _{j}$ is an eigenvalue of $A$.

(b) After an observation/measurement of an outcome $sh\left( \mu _{j}\right) 
$, the system is left in an eigenstate $\psi _{j}$ associated with the
detected eigenvalue $\mu _{j}$.

(b) In a measurement the transition probability $\mathcal{P}$ from a state $%
\psi $ to an eigenstate $\psi _{j}$ is given by 
\begin{equation*}
\mathcal{P}=\left\vert \left( \psi ,\psi _{j}\right) \right\vert ^{2}.
\end{equation*}

\textbf{Axiom U3}. The evolution of the state of a system is given by the
Shroedinger equation%
\begin{equation}
i\frac{\partial \psi }{\partial t}=H\psi  \label{sh}
\end{equation}%
where $H,$ the Hamiltonian operator, representing the energy of the system.

\bigskip

\textbf{Axiom U4}. Only the physical states represented by standard vectors
(namely vectors in $V$) can be produced in laboratory.

\subsection{Discussion of the axioms}

AXIOM 1. In the classical formalism, a physical system is not described only
by a given Hilbert space as axiom C1 claims, but by an Hilbert space and the
domain of a self-adjoint realization of the Hamiltonian operator. On the
contrary, in the ultravectors formalism the physical system is described
just by the space $V_{\mathcal{B}}$. Let see an example:

\textbf{A particle in a box. }For simplicity, we consider a one-dimensional
model and suppose that the box is modelled by the interval $\left[ 0,1\right]
.$ Clearly, the Hilbert space $L^{2}\left( 0,1\right) $ is not sufficient to
describe the system but it is necessary to give the Hamiltonian%
\begin{equation*}
H:H^{2}\left( 0,1\right) \cap H_{0}^{1}\left( 0,1\right) \rightarrow
L^{2}\left( 0,1\right)
\end{equation*}%
defined by 
\begin{equation}
H\psi =-\frac{1}{2m}\Delta \psi  \label{hhh}
\end{equation}%
where $\Delta \psi $ must be intended in the sense of distribution (here $m$
denotes the mass of the particle and we have assumed $\hslash =1$).

\textbf{A particle in a ring. }Now suppose that a point-particle is
constrained in a ring of length 1. Also in this case any state can be
represented by a vector in the Hilbert space $L^{2}\left( 0,1\right) ,$ but
in order to describe the system is necessary to give a different selfadjoint
realization of the Hamiltonian operator, namely an operator having the form (%
\ref{hhh}), but defined on the domain%
\begin{equation*}
H:H_{per}^{2}\left( 0,1\right) \rightarrow L^{2}\left( 0,1\right)
\end{equation*}%
where $H_{per}^{2}\left( 0,1\right) $ is the closure in the $H^{2}$ norm of
the space 
\begin{equation*}
\mathcal{C}_{per}^{2}\left[ 0,1\right] =\left\{ \psi \in \mathcal{C}^{2}(%
\left[ 0,1\right] ,\mathbb{C})\ |\ \psi (0)=\psi (1);\ \psi ^{\prime
}(0)=\psi ^{\prime }(1)\right\}
\end{equation*}

\bigskip

Now let us see how these two cases can be described in the ultrafunctions
formalism.

\bigskip

\textbf{A particle in a box. }In this case, the system is described by the
space%
\begin{equation*}
V_{\mathcal{B}}^{2,0}\left[ 0,1\right] :=\mathcal{B}\left[ \mathcal{C}%
_{0}^{2}\left[ 0,1\right] \right]
\end{equation*}%
The Hamiltonian operator $H$ is given by the canonical extension of $-\frac{1%
}{2m}\Delta $ to $\mathcal{B}\left[ \mathcal{C}_{0}^{2}\left[ 0,1\right] %
\right] $.

\textbf{A particle in a ring. }In this case, the system is described by the
space%
\begin{equation*}
V_{\mathcal{B}}^{2,per}\left[ 0,1\right] :=\mathcal{B}\left[ \mathcal{C}%
_{per}^{2}\left[ 0,1\right] \right]
\end{equation*}%
and the Hamiltonian operator $H$ is given by the canonical extension of $-%
\frac{1}{2m}\Delta $ to $\mathcal{B}\left[ \mathcal{C}_{per}^{2}\left[ 0,1%
\right] \right] $.

\bigskip

Thus in the ultrafunctions description, different physical systems give
different ultrafunction spaces; on the contrary, the Hamiltonian is given by
the unique canonical extension of $-\frac{1}{2m}\Delta $ in the relative
spaces.

\bigskip

AXIOM 2. In the ultrafunction formalism, the notion of self-adjoint operator
is not needed. In fact osservables can be represented by internal Hermitian
operators. It follows that any observable has exactly $\beta =\dim ^{\ast
}(V_{\mathcal{B}})$ eigenvalues (of course, if you take account of their
multiplicity). No essential distinction between eigenvalues and continuous
spectrum is required. For example, consider the eigenvalues of the position
operator $\hat{q}$ of a free particle. The eigenfunction relative to an
eigenvalue $q\in \mathbb{R}$ is an ultrafuncion of $\delta $-type
concentrated at the point $q$ (see Def. \ref{dt}).

In general the eigenvalues $\mu $'s of an internal Hermitian operator $A$
are hyperreal numbers, and hence, assuming that a measurement gives a real
number, we have imposed in Axiom 2 that the outcome of an experiment is $%
sh(\mu )$. However, we think that the probability is better described by the
hyperreal number $\left\vert \left( \psi ,\psi _{j}\right) \right\vert ^{2}$
rather than the real number $sh(\left\vert \left( \psi ,\psi _{j}\right)
\right\vert ^{2})$ (see \cite{BHW} for a presentation and discussion of the
Non Archimedean Probability). For example, let $\psi \in \mathcal{D}$ be the
state of a system; the probability of finding a particle in the position $q$
is given by%
\begin{equation*}
\left\vert \int \psi (x)\eta e_{q}(x)dx\right\vert =\eta \left\vert \psi
(q)\right\vert
\end{equation*}%
where $e_{q}$ is a $\delta $-type function and the normalization factor%
\begin{equation*}
\eta =\left\Vert e_{q}\right\Vert _{\left( L^{2}\right) ^{\ast }}^{-1}\sim 0
\end{equation*}%
is an infinitesimal number.

\bigskip

AXIOM 3. Since $H$ is an internal operator defined on a hyperfinite vector
space it can be represented by an Hermitian hyperfinite matrix and hence the
evolution operator of (\ref{sh}) is the exponential matrix $e^{tH}.$

\bigskip

AXIOM 4. In ultrafunction theory, the mathematical distinction between the
standard states and the ideal states is intrinsic and it does not correspond
to anything in the usual formalism. The point is to know if it corresponds
to something physically meaningful. Basically, we can say that the standard
states can be prepared in a laboratory, while the ideal states represent
"extreme" situations useful in the foundations of the theory and in thought
experiments (gedankenexperiment). For example the Dirac $\delta $-measure is
not a standard state but an ideal state and it represents a situation in
which the position of a particle is perfectly determined. Clearly this
situation cannot be produced in a laboratory, but nevertheless it is useful
in our description of the physical world. The standard states are
represented by functions in $V$ which is chosen depending on the model of
the physical system. The other states (namely, the states in $V_{\mathcal{B}%
}\backslash V$) are the ideal states. This situation makes more explicit
something which is already present in the classical approach. For example,
in the Shroedinger representation of a free particle in $\mathbb{R}^{3}$,
consider the state 
\begin{equation*}
\psi (x)=\frac{\varphi (x)}{|x|},\ \varphi \in \mathcal{D}(\mathbb{R}^{3}),\
\varphi (0)>0.
\end{equation*}%
We have that $\psi (x)\in L^{2}(\mathbb{R}^{3})$ but this state cannot be
produced in a laboratory, since the expected value of its energy 
\begin{equation*}
\left( H\psi ,\psi \right) =\frac{1}{2m}\int \left\vert \nabla \psi
\right\vert ^{2}dx
\end{equation*}%
is infinite. In other words, Axiom 4 makes formally precise something which
is already present (but hidden) in the classical theory. This point will be
discussed also in the next section.

\bigskip

\subsection{The Heisenberg algebra}

\bigskip

In this section we will apply ultrafunction theory to the description of a
quantum particle via the algebraic approach. For simplicity here we consider
the one-dimensional case. The states of a particle are defined by the
observables $q$ and $p$ which represent the position and the momentum
respectively. A quantum particle is described by the algebra of observables
generated by $p$ and $q$ according to the following commutation rules:%
\begin{equation*}
\left[ p,q\right] =i,\ \ \left[ p,p\right] =0,\ \ \left[ q,q\right] =0
\end{equation*}%
The algebra generated by $p$ and $q$ with the above relations is called the
Heisenberg algebra and denoted by $\mathfrak{A}_{H}$. The Heisenberg algebra
does not fit in the general theory of $C^{\ast }$-algebras since both $p$
and $q$ are not bounded operator. The usual technical solution to this
problem is done via the Weyl operators and the Weyl algebra (for more
details and a discussion on this point we refer to \cite{strocchi05}).

Let us see an alternative approach via ultrafunction theory. First of all we
take a representation of $\mathfrak{A}_{H},$ namely an algebra homomorphism%
\begin{equation*}
J:\mathfrak{A}_{H}\rightarrow \mathfrak{L}(V)
\end{equation*}%
where $\mathfrak{L}(V)$ is the algebra of the linear operators on a complex
vector space $V\subset H\in \mathcal{U}$ where $H$ is an Hilbert space and $%
\mathcal{U}$ is our universe (see section \ref{OL}). To fix the ideas, we
can consider the following "classical example":%
\begin{equation*}
H=L^{2}(\mathbb{R});\ \ \ V=\mathcal{S};
\end{equation*}%
\begin{equation*}
J(p)=-i\partial ;\ \ J(q)=x.
\end{equation*}

The quantum system of a particle will be described by the ultravector space $%
V_{\mathcal{B}}=\mathcal{B}\left[ V\right] $. The operators $J(p)$ and $J(q)$
can be extended to the space $V_{\mathcal{B}}$ according to definition (\ref%
{CE}); such extensions will be called $\hat{p}$ and $\hat{q}$ respectively. $%
\hat{p}$ and $\hat{q}$ are Hermitian operators and hence $V_{\mathcal{B}}$
has an othonormal basis generated by the eigenfunctions of $\hat{p}$ or $%
\hat{q}$. Let $\left\{ e_{a}\right\} _{a\in \Sigma }$ be the eigenfunctions
of $\hat{q}$ corresponding to the eigenvalue $a\in \Sigma \subset \mathbb{R}%
^{\ast }$. A very interesting fact is that the eigenfunctions violate the
Heisenberg relation $\left[ \hat{p},\hat{q}\right] =i.$

To see this fact we argue indirectly. Assume that the Heisenberg relation
holds; then 
\begin{equation*}
\left( \left[ \hat{p},\hat{q}\right] e_{a},e_{a}\right) =i\left\Vert
e_{a}\right\Vert ^{2}.
\end{equation*}

On the other hand, by a direct computation, we get:%
\begin{eqnarray*}
\left( \left[ \hat{p},\hat{q}\right] e_{a},e_{a}\right) &=&\left( \left( 
\hat{p}\hat{q}-\hat{q}\hat{p}\right) e_{a},e_{a}\right) =\left( \hat{p}\hat{q%
}e_{a},e_{a}\right) -\left( \hat{q}\hat{p}e_{a},e_{a}\right) \\
&=&\left( \hat{q}e_{a},\hat{p}e_{a}\right) -\left( \hat{p}e_{a},\hat{q}%
e_{a}\right) =a\left( e_{a},\hat{p}e_{a}\right) -a\left( \hat{p}%
e_{a},e_{a}\right) =0.
\end{eqnarray*}

This fact is consistent with the Axiom \textbf{U4 }which establishes that
the ideal states cannot be produced in laboratory. According to this
description of QM, the uncertainty relations hold only for the limitation of
the experimental apparatus. In a laboratory you can prepare a state
corresponding to a function $\psi $ in the space $V=\mathcal{S}$, but you
cannot prepare a state such as $e_{a}\in V_{\mathcal{B}}\backslash \mathcal{S%
}$ which corresponds to a particle which is exactly in the position $a.$

\bigskip

\end{document}